\documentclass[12pt,leqno]{article}
\usepackage{amsmath,amssymb,amsthm,amsfonts,mathrsfs,amscd}
\usepackage{latexsym,enumerate,color,geometry}
\usepackage{hyperref}

\def\d{\mathrm{d}}     

 \newtheorem{lemma}{Lemma}[section]
 \newtheorem{theorem}[lemma]{Theorem}
 \newtheorem{corollary}[lemma]{Corollary}
 \newtheorem{remark}[lemma]{Remark}

 \numberwithin{equation}{section}

\def\<{\langle}
\def\>{\rangle}
\def\E{\mathbb{E}}
\def\Dom{\mathscr{D}}
\def\Ran{\mathrm{Ran}}

\def\e{\operatorname{e}}

\def\de{\end{equation}}

\title{{\bf Harnack Inequality for Semilinear SPDE with Multiplicative Noise}\footnote{Supported by NSFC(11131003), SRFDP, 985-Project.}}

\author{
{\bf Zhang Shao-Qin}\\
\footnotesize{School of Math. Sci. and Lab. Math. Com. Sys., Beijing Normal University, Beijing 100875, China}\\
\footnotesize{Email: zhangsq@mail.bnu.edu.cn}
}

\begin{document}
\maketitle

\begin{abstract}
By a new approximate method, dimensional free Harnack inequalities are established for a class of semilinear stochastic differential equations in Hilbert space with multiplicative noise. These inequalities are applied to study the strong Feller property for the semigroup and some properties of invariant measure.
\end{abstract}\noindent

\noindent


\vskip 2cm

\section{Introduction and main results}
The main aim of this paper is to prove Harnack inequality for semilinear stochastic equations on Hilbert spaces with multiplicative noise. This type of inequality, which was proved for the first time in \cite{Wang97}, has became a powerful tool in infinite dimensional stochastic analysis. There are many papers prove this type of inequality for SPDE with additive noise, see \cite{DPRWang2009,Liu09,LiuW08,Ouyang2009a,Ouyang2009b,Ouyang2011a,OuyangRW2012,Wang2007,Wang2011,WangWX2011,WangX2011} and reference therein. In \cite{RoWang2010}, the log-Harnack inequality for semilinear SPDE with non-additive noise was proved for the first time, but by the gradient estimate method used there, only determine and time independent coefficient was treated. A new method to deal with the case of general coefficients for SDE was introduced in \cite{Wang2011}. This method has been generalized to functional stochastic differential equations, see \cite{WangY2011}. In this paper, we generalized this method to the case of semilinear SPDE. There are some disadvantages for finite dimension approximate method here, see Remark \ref{remark2}, therefore we use the coupling argument again as in \cite{Wang2011} with a slight modification. Since it seems not so clear to solves the similar equation of process $Y_{t}$ ( see equation (2.3) in \cite{Wang2011} ) in infinite dimension, we turn to a new process which plays the role as the difference of the coupling processes, we get it as a local strong solution of a SPDE and solve the equation by truncation in the same sprite in \cite{Brzez97}. By this process and Girsanov theorem, we get a coupling in a new probability space. On the other hand, we get Harnack inequality by another type of  approximation. We perturb the linear term by a suitable linear operator which closely relates to diffusion term. It's different from finite dimensional approximate and Yosida approximate, by this perturbation, we get a stronger linear term and it makes us to prove the inequality for the perturbed equation more easy.

Let $H$ be a separable Hilbert space with inner product $\<\cdot,\cdot\>$ and norm $\|\cdot\|$, consider the following stochastic differential equation on $H$:
\begin{equation}\label{equ1}
        \d x_{t} = -Ax_{t}\d t + F(t,x_{t})\d t+B(t,x_{t})\d W_{t}
\end{equation}
$W=W(t),t\geq0$ is a cylindrical Brownian motion on $H$ with covariance operator $I$ on filtered probability space $(\Omega, \mathcal{F}, \mathbb{P},(\mathcal{F}_{t})_{t\geq0})$, and the coefficients satisfy the following hypotheses:
\begin{enumerate}[({H}1)]

\item $A$ is a negative self adjoint operator with discrete spectrum:
\begin{equation}
0\leq\lambda_{1} \leq \lambda_{2} \leq \cdots \leq \lambda_{n} \rightarrow \infty,
\end{equation}
$\{\lambda_{n},n\in\mathbb{N}\}$ are the eigenvalues of $A$, and $\{e_{n}\}_{n=1}^{+\infty}$ are the corresponding eigenvectors, the compact $C_{0}$ semigroup generated by $-A$ denoted by $S(t)$.\label{itemH1}

\item $F:[0,\infty)\times\Omega\times H \rightarrow H$ and $B: [0,\infty)\times\Omega\times H \rightarrow L(H)$ are $\mathscr{P}_{\infty}\times \mathscr{B}(H)$ measurable, here $\mathscr{P}_{\infty}$ is predictable $\sigma$-algebra on $[0,\infty)\times\Omega$ and $L(H)$ is all the bounded operators on $H$, and there exists an increasing function $K_{1}:[0,+\infty)\rightarrow [0,\infty)$, such that
\begin{equation}
||F(t,x)-F(t,y)||+||B(t,x)-B(t,y)||_{HS}\leq K_{1}(t)||x-y||,
\end{equation}
for all $t\geq 0$, $x \in H$, $\mathbb{P}$-a.s, here $||\cdot||_{HS}$ denote the Hilber-Schmidt norm, and there exists $r>1$, such that for all $t>0$,
\begin{equation}
\mathbb{E}\left(\int_{0}^{t}\left|\left|F(s,0)\right|\right|\d s\right)^{r} < \infty,
\end{equation}
\begin{equation}
\sup_{u\in [0,t]}\int_{0}^{u}\left(\mathbb{E}\left|\left|S(u-s)B(s,0)\right|\right|_{HS}^{2r}\right)^{\frac{1}{r}}\d s < \infty,
\end{equation}
\label{itemH2}

\item There exist a decreasing function $\rho:[0,\infty)\rightarrow(0,\infty)$, and a bounded self adjoint operator $B_{0}$ satisfying that there exists $\{b_{n}>0|n\in\mathbb{N}\}$ such that $B_{0}e_{n}=b_{n}e_{n}$ and
\begin{equation}
B(t,x)B(t,x)^{*} \geq \rho(t)^{2}B_{0}^{2},\ \forall x\in H, t\geq 0,\ \mathbb{P}\mbox{-a.s.},
\end{equation}
\label{itemH3}

\item $\textrm{Ran}(B(t,x)-B(t,y))\subset \mathscr{D}(B_{0}^{-1})$ holds for all $(t,x) \in [0,\infty)\times H, \mathbb{P}$-a.s., and there exists an increasing function $K_{2}:[0,\infty)\rightarrow\mathbb{R}$ such that
    \begin{align*}
    2\<F(t,x)-F(t,y),B_{0}^{-2}(x-y)\>+&||B_{0}^{-1}(B(t,x)-B(t,y))||_{HS}^{2}\\
    &\leq K_{2}(t)||B_{0}^{-1}(x-y)||^{2}
    \end{align*}
holds for all $x,y \in \mathscr{D}(B_{0}^{-2})$ and all $t\geq 0$, $\mathbb{P}$-a.s.,\label{itemH4}

\item There exists an increasing function $K_{3}:[0,\infty)\rightarrow (0,\infty)$, such that $||(B(t,x)^{*}-B(t,y)^{*})B_{0}^{-2}(x-y)||\leq K_{3}(t)||x-y||_{H_{0}}$ holds for all $x,y\in H$, $t\geq 0$ and $x-y\in\Dom(B_{0}^{-1})$ almost surely.\label{itemH5}
\end{enumerate}

\begin{remark}\label{remark1}
\begin{enumerate}[(1)]
\item Under (H\ref{itemH1}), we can replace $\mathscr{D}(B_{0}^{-2})$ in (H\ref{itemH4}) by  $\bigcup_{n} H_{n}$, where $H_{n}= \mathrm{span}\{e_{1},\cdots,e_{n}\}$.
\item  (H\ref{itemH3}) equals to that $\Ran(B(t,x)) \supset \Ran B_{0}$ and $||B(t,x)^{-1}z||\leq \rho(t)^{-1}||B_{0}^{-1}z||$, for all $z \in \mathscr{D}(B_{0}^{-1})$, $t\geq 0$, $\mathbb{P}$\mbox{-a.s.},
\item (H\ref{itemH5}) will be used as a condition in addition to get Harnack inequality, and by (H\ref{itemH4}), $B_{0}^{-1}(B(t,x)-B(t,y))$ is an bounded operator, so in (H\ref{itemH5}) we only require $x-y\in \Dom(B_{0}^{-1})$.
\end{enumerate}
\end{remark}
For the proof of Remark \ref{remark1}, see Appendix. We state our main result of this paper
\begin{theorem}\label{mainthm}
If (H\ref{itemH1})-(H\ref{itemH4}) hold, then
\begin{equation}
P_{T}\log{f}(y)\leq\log{P_{T}f(x)}+\frac{K_{2}(T)||x-y||_{H_{0}}}{2(1-e^{K_{2}T})},\ \forall f\in \mathscr{B}_{b}(H),f\geq 1, x,y\in H,T>0.
\end{equation}
If, in addition, (H\ref{itemH5}) holds, then for $p>(1+\frac{K_{3}(T)}{\rho(T)})^{2}$, $\delta_{p,T} = K_{3}\vee \frac{\rho(T)}{2}(\sqrt{p}-1)$, the Harnack inequality
\begin{equation}\label{Harineq}
(P_{T}f(y))^{p}\leq (P_{T}f^{p}(x))
\exp{\left[\frac{K_{2}(T)\sqrt{p}(\sqrt{p}-1)||x-y||_{H_{0}}^{2}}{4\delta_{p,T}[(\sqrt{p}-1)\rho(T)-\delta_{p,T}](1-e^{K_{2}T})}\right]},
\end{equation}
holds for all $T>0$, $x,y\in H$ and $f\in\mathscr{B}_{b}^{+}(H)$, where $||x||^{2}_{H_{0}}=\sum_{n=0}^{+\infty}{b_{n}^{-1}}\<x,e_{n}\>^{2}$, $H_{0}=\Dom(B_{0}^{-1})$.
\end{theorem}

\begin{remark}\label{remark2}
One may use the finite dimension approximate method to get the Harnack inequalities, but here we mention that there are difficulties to overcome and it may not be better than the method used here. Let $\pi_{n}$ be the projection form $H$ to $H_{n}$, then get the following equation on $H_{n}$
\begin{equation}
\d x^{n}_{t} = -A_{n}x^{n}_{t}\d t + F_{n}(t,x^{n}_{t})\d t+B_{n}(t,x^{n}_{t})\d W^{n}_{t},
\end{equation}
where,
\begin{equation}
A_{n}=\pi_{n}A,\ F_{n}=\pi_{n}F|_{H_{n}},\ B_{n}=\pi_{n}B|_{H_{n}},\ W^{n}=\pi_{n}W,
\end{equation}
one may find that after projecting to lower dimension, an invertible operator may become degenerate, for example, an operator has the matrix form,
$
\left(
        \begin{array}{cc}
          0 & 1 \\
          1 & 0 \\
        \end{array}
      \right)
$
, under the orthonormal basis $\{e_{1},\ e_{2}\}$. It's easy to find that it's degenerate after projecting to the subspace generated by $e_{1}$. By (H\ref{itemH3}), one may replace $B$ by its symmetrization $\sqrt{BB^{*}}$, but constant may become worse in (H\ref{itemH2}) and (H\ref{itemH4}), see remark after theorem 1 in \cite{ArY81}, and it seems not easy to get similar estimate for $\sqrt{BB^{*}}$ as in (H\ref{itemH4}).
\end{remark}

\section{Proof of Theorem \ref{mainthm}}
Fixed a time $T>0$, we focus our discussion on the interval [0,T]. In order to prove the main theorem, we need some lemmas, and denote $K_{i}(T)$ by $K_{i}$, $i=1,2,3$, for for simplicity's sake.  The first lemma prove the existence and uniqueness of mild solution of the equation (\ref{equ1}), and give some estimates.
\begin{lemma}
Under the condition (\textrm{H}\ref{itemH1}) and (\textrm{H}\ref{itemH2}), equation (\ref{equ1}) has a pathwise unique mild solution and
\begin{equation}\label{solution_estimate_1}
\sup_{t\in[0,T]}\E{||x_{t}||^{r}}\leq C(r,T)(1+\E||x_{0}||^{r}).
\end{equation}
\end{lemma}
\noindent\emph{Proof.}
The existence part goes along the same lines as that of Theorem 7.4 in \cite{DPZ1992}, if we can prove that there exists $p\geq2$, such that
\begin{equation}
\sup_{t\in[0,T]}\E\left|\left|\int_{0}^{t}{e^{-(t-s)A}F(s,x_{s})\d s}\right|\right|^{p}<\infty,
\end{equation}
and
\begin{equation}
\sup_{t\in[0,T]}\E{\left|\left|\int_{0}^{t}{e^{-(t-s)A}B(s,x_{s})\d W_{s}}\right|\right|^{p}}<\infty
\end{equation}
for all $H$-valued predictable processes $x$ defined on $[0,T]$ satisfying
\begin{equation}
\sup_{t\in[0,T]}\E{||x_{t}||^{p}}<\infty.
\end{equation}
In fact, for $r$ in (H\ref{itemH2}),
\begin{equation*}
\begin{split}
&\sup_{t\in[0,T]}\E{\left|\left|\int_{0}^{t}{e^{-(t-s)A}B(s,x_{s})\d W_{s}}\right|\right|^{r}}\\
\leq &\sup_{t\in[0,T]}{\E{\left|\left|\int_{0}^{t}{e^{-(t-s)A}(B(s,x_{s})-B(s,0))\d W_{s}}\right|\right|^{r}}}+\sup_{t\in[0,T]}{\E{\left|\left|\int_{0}^{t}{e^{-(t-s)A}B(s,0)\d W_{s}}\right|\right|^{r}}}\\
\leq & C(r,T)(1+\E{||x_{t}||^{r}})+\left(\frac{r}{2}(r-1)\right)^{\frac{r}{2}}
\sup_{t\in[0,T]}\left(\int_{0}^{t}\left(\E||S(t-s)B(s,0)||_{HS}^{r}\right)^{\frac{2}{r}}\right)^{r}
<\infty.
\end{split}
\end{equation*}
$F$ is treated similarly, we omit it. Estimate (\ref{solution_estimate_1}) follows from Grownwall's lemma. For the uniqueness part. If $x_{t}^{1},x_{t}^{2}$ are mild solutions of equation (\ref{equ1}), then
\begin{equation}
\begin{split}
\E\sup_{u\in[0,t]}{||x_{u}^{1}-x_{u}^{2}||^{r}}\leq& 2^{r}T\E{\sup_{u\in[0,t]}{\int_{0}^{u}{||S(u-s)(F(s,x_{s}^{1})-F(s,x_{s}^{2}))||^{r}\d s}}}\\
&+ 2^{r}\E\sup_{u\in[0,t]}{||\int_{0}^{u}{S(u-s)(B(t,x^{1}_{s})-B(t,x^{2}_{s}))\d W_{s}}||^{r}}\\
\leq& 2^{r}T\int_{0}^{t}{\E{||x_{u}^{1}-x_{u}^{2}||^{r}}\d s}+C(r,T)\E\int_{0}^{t}{||x_{s}^{1}-x_{s}^{2}||^{r}\d s}\\
\leq& C(r,T)\int_{0}^{t}{\E\sup_{u\in[0,s]}{||x_{u}^{1}-x_{u}^{2}||^{r}}\d s},
\end{split}
\end{equation}
by the second inequality, $\E\sup_{u\in[0,t]}{||x_{u}^{1}-x_{u}^{2}||^{r}}<\infty$, then by Gronwall's lemma, $x_{t}^{1}=x_{t}^{2},\ \forall t\in[0,T]$, $\mathbb{P}$-a.s.

\qed

Denote $A_{\epsilon} = A + \epsilon B_{0}^{-2}$, $\Dom(A_{\epsilon}) = \Dom(A) \bigcap \Dom(B_{0}^{-2}) \subset \Dom(B_{0}^{-2})$, it is a self adjoint operator, the eigenvalues of $A_{\epsilon}$ are $\{\lambda_{n,\epsilon}:=\lambda_{n}+\epsilon b_{n}^{-2}\ | n\in \mathbb{N}\}$ and the eigenvectors remain $\{e_{n}|n\in\mathbb{N}\}$. In fact, one can define a self adjoint operator $\tilde{A}$ by
\begin{align}
\Dom(\tilde{A})&=\left\{x\in H \ |\  \sum_{n=0}^{+\infty}{(\lambda_{n}+\epsilon b_{n}^{-2})^{2}\<x,e_{n}\>^{2}}<+\infty \right\},\\ \tilde{A}x&=\sum_{n=0}^{+\infty}{(\lambda_{n}+\epsilon b_{n}^{-2})\<x,e_{n}\>e_{n}},
\end{align}
then by basic inequality and spectral decomposition of $A$ and $B_{0}^{-2}$, it is easy to see that $\tilde{A}=A_{\epsilon}$.
\begin{lemma}
For the mild solution of equation
\begin{equation}\label{equ2}
\d x_{t}^{\epsilon}= -(A+\epsilon B_{0}^{-2})x_{t}^{\epsilon}\d t + F(t,x_{t}^{\epsilon})\d t + B(t,x_{t}^{\epsilon})\d W_{t},\ x_{0}^{\epsilon}= x,
\end{equation}
we have
\begin{equation}
\lim_{\epsilon \rightarrow 0^{+}}\mathbb{E}||x_{t}-x_{t}^{\epsilon}||^{2}=0,\ \forall t\in [0,T].
\end{equation}
\end{lemma}
\noindent\emph{Proof.}
Since
\begin{align}
x_{t} &= e^{-tA}x + \int_{0}^{t}{e^{-(t-s)A}F(s,x_{s})\d s} + \int_{0}^{t}{e^{-(t-s)A}B(s,x_{s})\d W_{s}},\\
x^{\epsilon}_{t} &= e^{-t(A+\epsilon B_{0}^{2})}x + \int_{0}^{t}{e^{-(t-s)(A+\epsilon B_{0}^{2})}F(s,x^{\epsilon}_{s})\d s} + \int_{0}^{t}{e^{-(t-s)(A+\epsilon B_{0}^{2})}B(s,x^{\epsilon}_{s})\d W_{s}},
\end{align}
then
\begin{equation}
\begin{split}
||x_{t}-x_{t}^{\epsilon}||^{2}\leq& 3||(e^{-t\epsilon B_{0}^{-2}}-1)e^{(-tA)}x||^{2}\\
&+3||\int_{0}^{t}{(e^{-(t-s)A}F(s,x_{s})-e^{-(t-s)(A+\epsilon B_{0}^{-2})}F(s,x_{s}^{\epsilon}))\d s}||^{2}\\
&+3||\int_{0}^{t}{(e^{-(t-s)A}B(s,x_{s})-e^{-(t-s)(A+\epsilon B_{0}^{-2})}B(s,x_{s}^{\epsilon}))\d W_{s}}||^{2}\\
=:& I_{1}+I_{2}+I_{3}.
\end{split}
\end{equation}
It's clear that $\lim_{\epsilon \rightarrow 0^{+}}I_{1}=0$. For $I_{2}$, we have
\begin{equation}
\begin{split}
I_{2}&\leq 6T\int_{0}^{t}{||(e^{-(t-s)A}-e^{-(t-s)(A+\epsilon B^{-2}_{0})})F(s,x_{s})||^{2}\d s}\\
&+6T\int_{0}^{t}{||e^{-(t-s)(A+\epsilon B^{-2}_{0})}(F(s,x_{s})-F(s,x_{s}^{\epsilon}))||^{2}\d s} =: I_{2,1}+I_{2,2},
\end{split}
\end{equation}
Since
\begin{align}
&||(e^{-(t-s)A}-e^{-(t-s)(A+\epsilon B_{0}^{-2})})F(s,x_{s})||\leq C(1+||x_{s}||),\\
&\lim_{\epsilon \rightarrow 0^{+}}{||(e^{-(t-s)A}-e^{-(t-s)(A+\epsilon B_{0}^{-2})})F(s,x_{s})||}=0.
\end{align}
By domain convergence theorem $\lim_{\epsilon \rightarrow 0^{+}}{\E I_{2,1}}=0$. On the other hand,
\begin{equation}
\begin{split}
I_{2,2}&\leq 6T\int_{0}^{t}{||e^{-(t-s)(A+\epsilon B_{0}^{2})}(F(s,x_{s})-F(s,x_{s}^{\epsilon}))||^{2}\d s}\\
&\leq 6T\int_{0}^{t}{||F(s,x_{s})-F(s,x_{s}^{\epsilon})||^{2}\d s}\leq 6TK_{1}\int_{0}^{t}{||x_{s}-x_{s}^{\epsilon}||^{2}\d s}.
\end{split}
\end{equation}
For $I_{3}$,
\begin{equation}
\begin{split}
\mathbb{E}{I_{3}}&\leq 6\mathbb{E}{||\int_{0}^{t}{(e^{-(t-s)A}-e^{-(t-s)(A+\epsilon B_{0}^{-2})})B(s,x_{s})\d W_{s}}||^{2}}\\
&+6\mathbb{E}{||\int_{0}^{t}{e^{-(t-s)(A+\epsilon B_{0}^{-2})}(B(s,x_{s})-B(s,x_{s}^{\epsilon}))\d W_{s}}||^{2}}=I_{3,1}+I_{3,2},
\end{split}
\end{equation}
and
\begin{equation}
\begin{split}
\E I_{3,1} \leq &12T \E ||\int_{0}^{t}{(I-e^{-(t-s)\epsilon B_{0}^{-2}})(e^{-(t-s)A}B(s,0))\d W_{s}}||^{2}\\
&+ 12T \E ||\int_{0}^{t}{(e^{-(t-s)A}-e^{-(t-s)(A+\epsilon B_{0}^{-2})})(B(s,x_{s})-B(s,0))\d W_{s}}||^{2}\\
\leq &12T \E \int_{0}^{t}{||(I-e^{-(t-s)\epsilon B_{0}^{-2}})(e^{-(t-s)A}B(s,0))||_{HS}^{2}\d s}\\
&+12T \E \int_{0}^{t}{||(I-e^{-(t-s)\epsilon B_{0}^{-2}})(e^{-(t-s)A}(B(s,x_{s})-B(s,0)))||_{HS}^{2}\d s}\\
=:&I_{3,1,1}+I_{3,1,2},
\end{split}
\end{equation}
since
\begin{align}
||(I-e^{-(t-s)\epsilon B_{0}^{-2}})e^{-(t-s)A}B(s,0)||^{2}=\sum_{n=1}^{+\infty}||(e^{-(t-s)\epsilon B_{0}^{-2}}-I)e^{-(t-s)A}B(s,0)e_{n}||^{2}
\end{align}
and
\begin{align}
&\lim_{\epsilon \rightarrow 0}||(e^{-(t-s)\epsilon B_{0}^{-2}}-1)e^{-(t-s)A}B(s,0)e_{n}||=0\\
&||(e^{-(t-s)\epsilon B_{0}^{-2}}-I)e^{-(t-s)A}B(s,0)e_{n}||\leq||e^{-(t-s)A}B(s,0)e_{n}||\\
\end{align}
and by (H\ref{itemH2})
\begin{equation}
\E \int_{0}^{t}{\sum_{n=1}^{+\infty}||e^{-(t-s)A}B(s,0)e_{n}||^{2}\d s}=\E \int_{0}^{t}{||e^{-(t-s)A}B(s,0)||_{HS}^{2}\d s}< \infty.
\end{equation}
By dominate convergence theorem, $\lim_{\epsilon \rightarrow 0}I_{3,1,1}=0$, Note that $B(s,x_{s})-B(s,0)\in L_{HS}(H)$, and
\begin{equation}
\begin{split}
&||(I-e^{-(t-s)\epsilon B_{0}^{2}})e^{-(t-s)A}(B(s,x_{s})-B(s,0))||_{HS}^{2}\\
=&\sum_{n=1}^{+\infty}{||(I-e^{-(t-s)\epsilon B_{0}^{2}})e^{-(t-s)A}(B(s,x_{s})-B(s,0))e_{n}||^{2}}
\end{split}
\end{equation}
and
\begin{align}
||(I-e^{-(t-s)\epsilon B_{0}^{-2}})(e^{-(t-s)A}(B(s,x_{s})-B(s,0)))e_{n}||^{2}&\leq ||(B(s,x_{s})-B(s,0))e_{n}||^{2}\\
\E \int_{0}^{t}{\sum_{n=1}^{+\infty}||(B(s,x_{s})-B(s,0))e_{n}||^{2}\d s}&\leq \E\int_{0}^{t}{||x_{s}||^{2}\d s}<\infty,
\end{align}
by dominate convergence theorem, $\lim_{\epsilon \rightarrow 0}{\E I_{3,1}}=0$.
Finally,
\begin{equation}
\E I_{3,2} \leq 6T \E\int_{0}^{t}{||B(s,x_{s})-B(s,x_{s}^{\epsilon})||_{HS}^{2}\d s}\leq 6TK_{2}\E\int_{0}^{t}{||x_{s}-x_{s}^{\epsilon}||^{2}\d s}.
\end{equation}
Now, we have
\begin{equation}
\E||x_{t}-x_{t}^{\epsilon}||^{2}\leq \psi_{\epsilon}(t) + C(T,K_{2})\E\int_{0}^{t}{||x_{s}-x_{s}^{\epsilon}||^{2}\d s}
\end{equation}
for some $\psi_{\epsilon}(t)$, which satisfies $\lim_{\epsilon \rightarrow 0}{\psi_{\epsilon}(t)}=0$, then by Gronwall's lemma,
\begin{equation}
\lim_{\epsilon \rightarrow 0}{\E||x_{t}-x_{t}^{\epsilon}||^{2}}=0,\ \forall t\in [0,T].
\end{equation}
\qed

Firstly, we shall consider the following equation, $\xi_{t}=\frac{2-\theta}{K_{2}}(1-\e^{K_{2}(t-T)})$,
\begin{equation}\label{equ3}
\begin{split}
\d z_{t}=&-A_{\epsilon}z_{t}\d t+(F(t,x_{t})-F(t,x_{t}-z_{t}))\d t + (B(t,x_{t})-B(t,x_{t}-z_{t}))\d W_{t}\\
&-\frac{1}{\xi_{t}}(B(t,x_{t}-z_{t})-B(t,x_{t}))B(t,x_{t})^{-1}z_{t}\d t - \frac{1}{\xi_{t}}z_{t}\d t,\  z_{0}=z.
\end{split}
\end{equation}
Note that, by (H\ref{itemH2})--(H\ref{itemH4}),
\begin{align}
&F(t,x_{t}) - F(t,x_{t}-z_{t}) \in H,\ (B(t,x_{t})-B(t,x_{t}-z_{t}))\in L_{HS}(H,H_{0}),\\
&(B(t,x_{t}-z_{t})-B(t,x_{t}))B(t,x_{t})^{-1}\in L(H_{0},H_{0}),
\end{align}
it's natural to solve the equation in $H_{0}$, we shall search a suitable Gelfand triple. To this end, we should restrict the operator $A_{\epsilon}$ to $H_{0}$.

\begin{lemma}
Define $A_{0,\epsilon}$ as follows
\begin{align}
\Dom(A_{0,\epsilon}) = B_{0}(\Dom(A_{\epsilon})),\ A_{0,\epsilon}x=A_{\epsilon}x, \forall x\in B_{0}(\Dom(A_{\epsilon})),
\end{align}
then, $A_{0,\epsilon}$ is well defined and $(A_{0,\epsilon}, B_{0}(\Dom(A_{\epsilon})))= (B_{0}A_{\epsilon}B_{0}^{-1},B_{0}(\Dom(A_{\epsilon})))$.
\end{lemma}
\noindent\emph{Proof.}
It's well defined. In fact for all $x\in B_{0}(\Dom(A_{\epsilon}))$,
\begin{equation}
\sum_{n=1}^{+\infty}{\lambda_{n,\epsilon}\<x,e_{n}\>^{2}}
=\sum_{n=0}^{+\infty}{\lambda^{2}_{n,\epsilon}b_{n}^{2}\<B_{0}^{-1}x,e_{n}\>^{2}}
\leq ||B||_{H}^{2}\sum_{n=1}^{+\infty}(\lambda_{n,\epsilon}^{2})\<B_{0}^{-1}x,e_{n}\>^{2}<+\infty,
\end{equation}
then $x\in\Dom(A_{\epsilon})$, and
\begin{align}
\sum_{n=1}^{+\infty}b_{n}^{-2}\<A_{\epsilon}x,e_{n}\>^{2}=\sum_{n=1}^{+\infty}{\lambda^{2}_{n,\epsilon}\<B_{0}^{-1}x,e_{n}\>^{2}}<+\infty,
\end{align}
then $A_{\epsilon}x\in\Dom(B_{0}^{-1}),\ \forall x\in B_{0}(\Dom(A_{\epsilon}))$, i.e. $A_{\epsilon}x\in H_{0}$. Finally, for all $x\in B_{0}(\Dom(A))$,
\begin{equation}
B_{0}A_{\epsilon}B_{0}^{-1}x=A_{\epsilon}B_{0}B_{0}^{-1}x=A_{\epsilon}x=A_{0,\epsilon}x.
\end{equation}
\qed

Now, we can define our Gelfand triple. Let
\begin{equation}
(V,||\cdot||_{V} )= (\Dom(A_{0,\epsilon}^{\frac{1}{2}}),||A_{0,\epsilon}^{\frac{1}{2}}\cdot||_{H_{0}}),
\end{equation}
then $(V^{*},||\cdot||_{V^{*}})$ is the complete of $(H_{0}, ||A_{0,\epsilon}^{-\frac{1}{2}}\cdot||_{H_{0}})$, $V^{*}\supset H_{0}\supset V$ is the triple we need. Since $\Dom(A_{\epsilon})\subset\Dom(B_{0}^{-2})$, $\Dom(A_{0,\epsilon})\subset \Dom(B_{0}^{-3})$, we have the following relationship moreover
\begin{equation}
V^{*}\supset H \supset H_{0} \supset \Dom(B_{0}^{-2}) \supset V.
\end{equation}

\begin{lemma}\label{lemma_strongsolution}
If conditions (H\ref{itemH1})-(H\ref{itemH4}) hold, equation (\ref{equ3}) has a unique strong solution up to the explosion time $\tau$.
\end{lemma}
\noindent\emph{Proof.}
Let
\begin{equation}
G_{n}(t,v)= \left\{
\begin{array}{ll}
    B(t,x_{t})^{-1}v,& \ ||v||_{H_{0}}\leq n,\\
    B(t,x_{t})^{-1}\frac{nv}{\ ||v||_{H_{0}}},& \ ||v||_{H_{0}}>n,\\
\end{array}
\right.
\end{equation}
and for simplicity's sake, we denote
$$F(t,x_{t}-v_{1})-F(t,x_{t}-v_{2}),\ G_{n}(t,v_{1})-G_{n}(t,v_{2}),\ B(t,x_{t})-B(t,x_{t}-z_{t})$$
by $F(t,v_{2},v_{1})$, $G_{n}(t,v_{1},v_{2})$, $\hat{B}(t,z_{t})$ respectively.
We consider the following equation firstly,
\begin{equation}\label{equ4}
\begin{split}
\d z_{t}=&-A_{0,\epsilon}z_{t}\d t + F(t,z_{t},0)\d t- \frac{1}{\xi_{t}}z_{t}\d t+\frac{1}{\xi_{t}}\hat{B}(t,z_{t})G_{n}(t,z_{t})\d t + \hat{B}(t,z_{t})\d W_{t}\\
=:&A_{n,\epsilon}(t,z_{t})\d t + \hat{B}(t,z_{t})\d W_{t}
\end{split}
\end{equation}
It's clearly that the hemicontinuous  holds, since $G_{n}(t,\cdot)$ remains a Lipschitz mapping from $H_{0}$ to $H$. By the direct calculus, see Appendix, we get that, for all $v,v_{1},v_{2} \in V$,
\begin{enumerate}[({A}1)]
\item Local monotonicity
\begin{equation*}
\begin{split}
&2 _{V^{*}}\<A_{n,\epsilon}(t,v_{1})-A_{n,\epsilon}(t,v_{2}),v_{1}-v_{2}\>_{V}
+||\hat{B}(t,v_{2})-\hat{B}(t,v_{1})||_{L_{HS}(H,H_{0})}^{2}\\
\leq&\left[K_{2}+\frac{2n\sqrt{K_{2}}-2}{\xi_{t}}
+\frac{n^{2}K_{1}||B_{0}||^{2}}{\epsilon^{2}\xi_{t}^{2}\delta^{2}}
+\frac{2}{\xi_{t}}(\sqrt{K_{2}}||v_{2}||_{H_{0}}^{2}
+\sqrt{\frac{2K_{1}}{\epsilon}}||B_{0}||\cdot||v_{2}||_{V}^{2})\right]\times\\
&\times||v_{1}-v_{2}||_{H_{0}}^{2}-2(1-\delta^{2})||v_{1}-v_{2}||_{V}^{2},\ \forall \delta\in(0,1).
\end{split}
\end{equation*}\label{itemA1}

\item Coercivity
\begin{equation*}
\begin{split}
&2 _{V^{*}}\<A_{n,\epsilon}(t,v),v\>_{V} + ||\hat{B}(t,v)||_{L_{HS}(H,H_{0})}^{2}\\
\leq& -2(1-\delta^{2})||v||_{V}^{2}+(\frac{n\sqrt{K_{2}}-2}{\xi_{t}}+\frac{n^{2}K_{1}}{\epsilon^{2}\xi_{t}^{2}\delta^{2}})||v||_{H_{0}}^{2}, \forall \delta\in(0,1).
\end{split}
\end{equation*}\label{itemA2}

\item Growth
\begin{equation*}
||A_{n,\epsilon}(t,v)||_{V^{*}}^{2}\leq \left(\frac{||B_{0}||^{2}}{\epsilon\xi_{t}}K_{2}
+\left(1+\frac{||B_{0}||^{4}K_{1}}{\epsilon\xi_{t}^{2}}\right)||v||_{V}^{2}\right)(1+||v||_{H_{0}}^{4}).
\end{equation*}\label{itemA3}
\end{enumerate}
Since
\begin{equation}\label{inequ1}
\begin{split}
||\hat{B}(t,v)||^{2}_{L_{HS}}=||B_{0}^{-1}\hat{B}(t,v)||^{2}_{HS}
\leq K_{2}||v||_{H_{0}}^{2}+\frac{2K_{1}}{\epsilon}||B_{0}||^{3}||v||_{V}||v||_{H_{0}}.
\end{split}
\end{equation}
does not satisfies the condition (1.2) in \cite{LiuR10}, but by the basic inequality
one can check that the proof in Lemma2.2 goes on well, see Appendix B. By the estimates above and Theorem 1.1 in \cite{LiuR10} for any $T_{0}<T$, equation (\ref{equ4}) has unique strong solution $(z_{t}^{n})_{t\in[0,T_{0}]}$, one can extends the solution to the interval $[0,T)$ by the pathwise uniqueness and continuous. Next we shall let $n$ goes to infinite. Let, $m>n$,
\begin{equation}
\tau_{m}^{n} = \inf\{t\in[0,T)\ |\ ||z_{t}^{m}||_{H_{0}}>n\},
\end{equation}
definite $\inf\emptyset=T$, then
\begin{equation}
\begin{split}
z_{t}^{m}=&z_{0} + \int_{0}^{t}{(-A_{0,\epsilon}z_{s}^{m} + F(s,z_{s}^{m},0)-\frac{1}{\xi_{s}}z_{s}^{m})\d s}\\
&-\int_{0}^{t}{\frac{1}{\xi_{s}}\hat{B}(s,z_{s}^{m})B(s,x_{s})^{-1}z_{s}^{m}\d s}
+\int_{0}^{t}{\hat{B}(s,z_{s}^{m})\d W_{s}},\ t<\tau_{m}^{n},
\end{split}
\end{equation}
by It\^{o}'s formula and (A\ref{itemA1}), for $t<\tau_{n}^{n}\wedge\tau_{m}^{n}$, we have
\begin{equation*}
\begin{split}
&\d ||z_{t}^{n}-z_{t}^{m}||_{H_{0}}^{2}-2\<\hat{B}(t,z_{t}^{n})-\hat{B}(t,z_{t}^{m}))\d W_{t}, z_{t}^{n}-z_{t}^{m}\>_{H_{0}}\\
=&\ 2 _{V^{*}}\<A_{n,\epsilon}(t,z_{t}^{n})-A_{n,\epsilon}(t,z_{t}^{m}),z_{t}^{n}-z_{t}^{m}\>_{V}
+||\hat{B}(t,z_{t}^{n})-\hat{B}(t,z_{t}^{m})||_{L_{HS}(H,H_{0})}\d t\\
\leq &\left(K_{2}+\frac{2}{\xi_{t}}(n\sqrt{K_{1}}+\sqrt{K_{2}}||z_{t}^{n}||_{H_{0}}^{2}
+\sqrt{\frac{2K_{1}}{\epsilon}}||B_{0}||\cdot||z_{t}^{n}||_{V}^{2})
+\frac{n^{2}K_{1}}{\epsilon^{2}\xi_{t}^{2}\delta^{2}}||B_{0}||^{2}\right)||z_{t}^{n}-z_{t}^{m}||_{H_{0}}^{2}\\
\end{split}
\end{equation*}
define
\begin{equation}
\begin{split}
&\Psi_{s}=K_{2}+\frac{2}{\xi_{s}}(\sqrt{K_{2}}||z_{s}^{n}||_{H_{0}}^{2}+n\sqrt{K_{1}}
+\sqrt{\frac{2K_{1}}{\epsilon}}||B_{0}||^{2}||z_{s}^{n}||_{V}^{2})
+\frac{n^{2}K_{1}||B_{0}||^{2}}{\epsilon^{2}\xi_{s}^{2}\delta^{2}},
\end{split}
\end{equation}
then
\begin{equation}
\begin{split}
&\exp{\left[-\int_{0}^{t}{\Psi_{s}\d s}\right]}||z_{t}^{n}-z_{t}^{m}||_{H_{0}}^{2}\\
\leq &\int_{0}^{t}{2\exp{\left[-\int_{0}^{r}{\Psi_{s}\d s}\right]}\<(\hat{B}(r,z_{r}^{n})-\hat{B}(t,z_{r}^{m}))\d W_{r},z_{r}^{n}-z_{r}^{m}\>_{H_{0}}},
\end{split}
\end{equation}
therefore
\begin{equation}
\E\left\{\exp{\left[-\int_{0}^{t\wedge\tau_{n}^{n}\wedge\tau_{m}^{n}}{\Psi_{s}\d s}\right]}||z_{t\wedge\tau_{n}^{n}\wedge\tau_{m}^{n}}^{n}-z_{t\wedge\tau_{n}^{n}\wedge\tau_{m}^{n}}^{m}||_{H_{0}}^{2}\right\}=0.
\end{equation}
Note that
\begin{equation}
\E\int_{0}^{t}{||z_{s}^{n}||^{2}_{V}\d s}<\infty,\ \forall t<T
\end{equation}
implies
\begin{equation}
\int_{0}^{t}{||z_{s}^{n}||_{V}^{2}\d s}<\infty,\ \forall t\in[0,T),\ \mathbb{P}\mbox{-a.s.},
\end{equation}
then
\begin{equation}
z_{t\wedge\tau_{n}^{n}\wedge\tau_{m}^{n}}^{n}=z_{t\wedge\tau_{n}^{n}\wedge\tau_{m}^{n}}^{m},\ \forall t\in[0,T),\ \mathbb{P}\mbox{-a.s.},
\end{equation}
let $t\uparrow T$, by the continuity, we have
\begin{equation}
z_{\tau_{n}^{n}\wedge\tau_{m}^{n}}^{n}=z_{\tau_{n}^{n}\wedge\tau_{m}^{n}}^{m},\ \mathbb{P}\mbox{-a.s.}
\end{equation}
If $\tau_{n}^{n}<\tau_{m}^{n}$, $z_{\tau_{n}^{n}}^{n}=z_{\tau_{n}^{n}}^{m}\in \partial B_{n}^{H_{0}}(0)$, by the definition of $\tau_{m}^{n}$, it's a contradictory. Thus $\tau_{n}^{n}\geq\tau_{m}^{n}$, similarly, $\tau_{n}^{n}\leq\tau_{m}^{n}$, so $\tau_{n}^{n}=\tau_{m}^{n}$, $\mathbb{P}$-a.s. and $z_{\tau_{n}^{n}}^{n}=z_{\tau_{m}^{n}}^{m}$. Therefore, we can definite
\begin{equation}
z_{t}=z_{t}^{n},\ t<\tau_{n}^{n};\ \tau=\sup_{n}{\tau_{n}^{n}},
\end{equation}
$(z,\tau)$ is a strong solution of equation (\ref{equ3}). By the same method, we can prove the uniqueness easily.
\qed
\bigskip
\\\emph{Proof of Theorem \ref{mainthm}}. Let
\begin{align*}
\d \tilde{W}_{s} &= \d W_{s} + \frac{1}{\xi_{s}}B(s,x_{s})^{-1}z_{s}\d s,\ s<T\wedge\tau\\
R_{s}&=\exp{\left[-\int_{0}^{s}\xi_{t}^{-1}\<B(t,x_{t})^{-1}z_{t},\d W_{t}\>-\frac{1}{2}\int_{0}^{s}{\frac{||B(t,x_{t})^{-1}z_{t}||^{2}}{\xi_{t}}\d t}\right]},\ s<T\wedge\tau,\\
\tau_{n}&=\inf\{t\in[0,T)\ |\ ||z_{t}||_{H_{0}}>n\},\ \mathbb{Q}:= R_{T\wedge\tau}\mathbb{P},
\end{align*}
write the equation of $z$ in the form of $\tilde{W}$:
\begin{equation}\label{equ5}
\d z_{t} = -A_{0,\epsilon}z_{t}\d t+F(t,z_{t},0)\d t + \hat{B}(t,z_{t})\d \tilde{W}_{t}- \frac{1}{\xi_{t}}z_{t}\d t,
\end{equation}
By It'\^{o}'s formula and (H\ref{itemH4}), for $s\in[0,T)$, and for $t<\tau_{n}\wedge s$,
\begin{equation}
\begin{split}
\d ||z_{t}||_{H_{0}}^{2} = &-2||z_{t}||_{V}^{2}\d t + 2 _{V^{*}}\<F(t,z_{t},0),z_{t}\>_{V}\d t-\frac{2||z_{t}||_{H_{0}}^{2}}{\xi_{t}}\d t\\
&+||\hat{B}(t,z_{t})||_{L_{HS}(H,H_{0})}^{2}\d t+2\<\hat{B}(t,z_{t})\d \tilde{W},z_{t}\>_{H_{0}}\\
\leq&\ 2 \<F(t,z_{t},0),B_{0}^{-2}z_{t}\>\d t + ||\hat{B}(t,z_{t})||_{L_{HS}(H,H_{0})}^{2}\d t \\
&-\frac{2||z_{t}||_{H_{0}}^{2}}{\xi_{t}}\d t +2\<\hat{B}(t,z_{t})\d \tilde{W},z_{t}\>_{H_{0}} \\
\leq& -\frac{2||z_{t}||_{H_{0}}^{2}}{\xi_{t}}\d t +2\<\hat{B}(t,z_{t})\d \tilde{W},z_{t}\>_{H_{0}}+ K_{2}||z_{t}||_{H_{0}}^{2}\d t,
\end{split}
\end{equation}
and
\begin{equation}
\begin{split}
\d \frac{||z_{t}||_{H_{0}}^{2}}{\xi_{t}} \leq & -\frac{2||z_{t}||_{H_{0}}^{2}}{\xi_{t}^{2}}\d t + \frac{K_{2}}{\xi_{t}}||z_{t}||_{H_{0}}^{2}\d t - \frac{\xi_{t}^{'}}{\xi_{t}^{2}}||z_{t}||_{H_{0}}^{2}\d t+\frac{2}{\xi_{t}}\<\hat{B}(t,z_{t})\d \tilde{W},z_{t}\>_{H_{0}}\\
=&\ \frac{2-K_{2}\xi_{t}+\xi_{t}^{'}}{\xi_{t}^{2}}||z_{t}||_{H_{0}}^{2}\d t+ \frac{2}{\xi_{t}}\<\hat{B}(t,z_{t})\d \tilde{W},z_{t}\>_{H_{0}}\\
=&\ \frac{\theta}{\xi_{t}^{2}}||z_{t}||_{H_{0}}^{2}\d t+ \frac{2}{\xi_{t}}\<\hat{B}(t,z_{t})\d \tilde{W},z_{t}\>_{H_{0}},
\end{split}
\end{equation}
by Girsanov theorem, $(\tilde{W})_{t<s\wedge\tau_{n}}$ is a Wiener process under the probability $\mathbb{Q}_{s,n} := R_{s\wedge\tau_{n}}\mathbb{P}$, and
\begin{equation}\label{inequ2}
\int_{0}^{s\wedge\tau_{n}}{\frac{||z_{t}||^{2}}{\xi_{t}^{2}}\d t}\leq \frac{||z_{0}||_{H_{0}}^{2}}{\theta\xi_{0}}+\int_{0}^{s\wedge\tau_{n}}{\frac{2}{\theta\xi_{t}}\<\hat{B}(t,z_{t})\d \tilde{W},z_{t}\>_{H_{0}}},
\end{equation}
then
\begin{equation}
\E_{\mathbb{Q}_{s,n}}{\int_{0}^{s\wedge\tau_{n}}{\frac{||z_{t}||^{2}}{\xi_{t}^{2}}\d t}}\leq \frac{||z_{0}||_{H_{0}}^{2}}{\theta\xi_{0}},
\end{equation}
Since, by (H\ref{itemH3})
\begin{equation}
\begin{split}
\log{R_{u}}&= -\int_{0}^{u}\xi_{t}^{-1}\<B(t,x_{t})^{-1}z_{t},\d \tilde{W}_{t}\>+\frac{1}{2}\int_{0}^{u}{\frac{||B(t,x_{t})^{-1}z_{t}||^{2}}{\xi_{t}}\d t}\\
&\leq -\int_{0}^{u}\xi_{t}^{-1}\<B(t,x_{t})^{-1}z_{t},\d \tilde{W}_{t}\>+\frac{1}{2\rho(T)^{2}}\int_{0}^{u}{\frac{||z_{t}||_{H_{0}}^{2}}{\xi_{t}}\d t},\ u\leq s\wedge\tau_{n},
\end{split}
\end{equation}

\begin{equation}
\E{R_{s\wedge\tau_{n}}}\log{R_{s\wedge\tau_{n}}}\leq \frac{||z_{0}||_{H_{0}}^{2}}{2\theta\xi_{0}\rho(T)^{2}},\ \forall s\in[0,T),n\geq1.
\end{equation}
As in \cite{Wang2011}, we can prove that $\{R_{s\wedge\tau}\ |\ s\in[0,T]\}$ is a martingale. Since
\begin{equation}
\E_{\mathbb{Q}}1_{[\tau_{n}\leq t]}\frac{||z_{t\wedge\tau_{n}}||_{H_{0}}^{2}}{\xi_{t\wedge\tau_{n}}}
\leq\E_{\mathbb{Q}}\frac{||z_{t\wedge\tau_{n}}||_{H_{0}}^{2}}{{\xi_{t\wedge\tau_{n}}}}
\leq\frac{||z_{0}||^{2}_{H_{0}}}{\xi_{0}},
\end{equation}
and
\begin{equation}
\E_{\mathbb{Q}}1_{[\tau_{n}\leq t]}\frac{||z_{t\wedge\tau_{n}}||_{H_{0}}^{2}}{\xi_{t\wedge\tau_{n}}}\geq\frac{n\mathbb{Q}(\tau_{n}\leq t)}{\xi_{0}}
\end{equation}
let $n$ goes to infinite, we have $\mathbb{Q}(\tau_{n}\leq t)=0$, $\forall t\in [0,T)$, then $\mathbb{Q}(\tau=T)=1$. Now, since $\tau = T$, $\mathbb{Q}$-a.s., equation (\ref{equ5}) can be solved up to time $T$. Let
\begin{equation}
\zeta=\inf\{t\in[0,T]\ |\ ||z_{t}||_{H_{0}}=0\},
\end{equation}
we shall prove that $\zeta \leq T$, here we assume $\inf\emptyset = +\infty$.
Otherwise, there exists a set $\Omega_{0}$, such that $\mathbb{P}(\Omega_{0})>0$, and for any $\omega \in \Omega_{0}$, $\zeta(\omega)>T$, then by the continuity of path, we have
\begin{equation}
\inf_{t\in[0,T]}{||z_{t}(\omega)||_{H_{0}}}>0,
\end{equation}
so
\begin{equation}
\int_{0}^{T}{\frac{||z_{t}||^{2}_{H_{0}}}{\xi_{t}^{2}}\d t}=+\infty,
\end{equation}
but
\begin{equation}
\E_{\mathbb{Q}}\int_{0}^{T}{\frac{||z_{t}||^{2}_{H_{0}}}{\xi_{t}^{2}}\d t}\leq \frac{||z_{0}||^{2}_{H_{0}}}{2\rho(T)^{2}\theta\xi_{0}}<+\infty,
\end{equation}
hence, $\zeta\leq T$, $\mathbb{Q}$-a.s., by the uniqueness of solution of equation (\ref{equ5}), we have
\begin{equation}
z_{t}\equiv 0,\ t>\zeta,\ \mathbb{Q}\mbox{-a.s.}
\end{equation}
Thus, $z_{T}=0$, $\mathbb{Q}$-a.s.

Next, we shall construct the coupling. Since under the probability space $(\Omega,\mathscr{F},R_{\tau\wedge T}\mathbb{P})$, $(\tilde{W}_{t})_{t\in[0,T]}$ is a Wiener process, let $y$ be the unique mild solution of the following equation
\begin{equation}
\d y_{t}=-A_{\epsilon}y_{t}\d t + F(t,y_{t})\d t + B(t,y_{t})\d \tilde{W}_{t},\ y_{0}=y,
\end{equation}
for $x_{t}$, it's the unique solution of the following equation
\begin{equation}
\d x_{t}=-A_{\epsilon}x_{t}\d t + F(t,x_{t})\d t -\frac{z_{t}}{\xi_{t}}\d t+ B(t,x_{t})\d \tilde{W}_{t},\ x_{0}=x.
\end{equation}
For the process $x_{t}-y_{t}$, it's the mild solution of the following equation
\begin{equation}\label{equ6}
\d u_{t}=-A_{\epsilon}u_{t}\d t+F(t,u_{t},0)\d t+\hat{B}(t,u_{t})\d \tilde{W}_{t} -\frac{z_{t}}{\xi_{t}}\d t,
\end{equation}
note that $z_{t}$ is a solution of equation
\begin{equation}
\d z_{t}=-A_{0,\epsilon}z_{t}\d t+F(t,z_{t},0)\d t+\hat{B}(t,z_{t})\d \tilde{W}_{t} -\frac{z_{t}}{\xi_{t}}\d t,
\end{equation}
Similar to equation (1.41), one can prove that equation (\ref{equ6}) has a strong solution in $H_{0}$,  since $V^{*} \supset H\supset H_{0}$ and $A_{0,\epsilon}$ is the restriction of $A_{\epsilon}$ to $H_{0}$, by the relation ship of variational solution and mild solution and the pathwise uniqueness, then $z_{t}=x_{t}-y_{t},\ \forall t\in[0,T]$, $\mathbb{Q}$-a.s.

By the method used in \cite{Wang2011}, we have log-Harnack inequality for equation (\ref{equ2}) :
\begin{equation}
\begin{split}
P_{T}^{\epsilon}\log{f}(y)&=\E_{\mathbb{Q}}\log{f(y_{T}^{\epsilon})}=\E{R_{T\wedge\tau}\log{f(x_{T}^{\epsilon})}}\leq\E{R_{T\wedge\tau}\log{R_{T\wedge\tau}}}+\log{\E{f(x_{T}^{\epsilon})}}\\
&\leq\log{P_{T}^{\epsilon}f(x)}
+\frac{||x-y||_{H_{0}}}{2\rho(T)^{2}\theta\xi_{0}}=\log{P_{T}^{\epsilon}f(x)}+\frac{K_{2}||x-y||_{H_{0}}}{2\rho(T)^{2}\theta(2-\theta)(1-e^{K_{2}T})}\ ,
\end{split}
\end{equation}
then by lemma 1.2, let $\epsilon\rightarrow 0$, and choose $\theta = 1$, for $f\in\mathscr{B}_{b}^{+}(H)$ and $f\geq 1$,
\begin{equation}
P_{T}\log{f}(y)\leq\log{P_{T}f(x)}+\frac{K_{2}||x-y||_{H_{0}}}{2\rho(T)^{2}(1-e^{K_{2}T})}.
\end{equation}
If (H\ref{itemH5}) holds in addition, by inequality (\ref{inequ2}), we have
\begin{equation}
\begin{split}
&\E_{s,n}{\exp{\left[h\int_{0}^{s\wedge\tau_{n}}{\frac{||z_{t}||^{2}_{H_{0}}}{\xi_{t}^{2}}\d t}\right]}}\\
\leq& \exp{\left[\frac{h||x-y||_{H_{0}}^{2}}{\theta\xi_{0}}\right]}
\E_{s,n}{\exp{\left[\frac{2h}{\theta}\int_{0}^{s\wedge\tau_{n}}{\frac{1}{\xi_{t}}\<\hat{B}(t,z_{t})\d \tilde{W},z_{t}\>}\right]}}\\
\leq&\exp{\left[\frac{h||x-y||_{H_{0}}^{2}}{\theta\xi_{0}}\right]}
\E_{s,n}{\left(\exp{\left[\frac{8h^{2}K_{3}^{2}}{\theta^{2}}
\int_{0}^{s\wedge\tau_{n}}{\frac{||z_{t}||^{2}_{H_{0}}}{\xi_{t}^{2}}\d t}\right]}\right)^{\frac{1}{2}}},
\end{split}
\end{equation}
for $h=\frac{\theta^{2}}{8K_{3}^{2}}$, and
\begin{equation}
\E_{s,n}{\exp{\left[\frac{\theta^{2}}{8K_{3}^{2}}\int_{0}^{s\wedge\tau_{n}}{\frac{||z_{t}||^{2}_{H_{0}}}{\xi_{t}^{2}}\d t}\right]}}
\leq \exp{\left[\frac{\theta K_{2}||x-y||_{H_{0}}^{2}}{4K_{3}^{2}(2-\theta)(1-e^{-K_{2}T})}\right]},
\end{equation}
Similar to \cite{Wang2011}, we get that
\begin{equation}
\sup_{s\in[0,T]}{\E{R_{s\wedge\tau}^{1+r}}\leq \exp{\left[\frac{\theta K_{2}(2K_{3}+\theta\rho(T))||x-y||_{H_{0}}}{8K_{3}^{2}(2-\theta)(K_{3}+\theta\rho(T))(1-e^{-K_{2}T})}\right]}}
\end{equation}
and for $p>(1+K_{3})^{2}$, $\delta_{p,T} = K_{3}\vee \frac{\rho(T)}{2}(\sqrt{p}-1)$,$f\in\mathscr{B}_{b}^{+}(H)$, choose $\theta = \frac{2K_{3}\rho(T)}{\sqrt{p}-1}$,
\begin{equation}
(P_{T}^{\epsilon}f(y))^{p}\leq (P_{T}^{\epsilon}f^{p}(x))
\exp{\left[\frac{K_{2}(T)\sqrt{p}(\sqrt{p}-1)||x-y||_{H_{0}}^{2}}{4\delta_{p,T}[(\sqrt{p}-1)\rho(T)-\delta_{p,T}](1-e^{K_{2}T})}\right]},
\end{equation}
by lemma 1.2, let $\epsilon\downarrow 0$, we have
\begin{equation}
(P_{T}f(y))^{p}\leq (P_{T}f^{p}(x))
\exp{\left[\frac{K_{2}(T)\sqrt{p}(\sqrt{p}-1)||x-y||_{H_{0}}^{2}}{4\delta_{p,T}[(\sqrt{p}-1)\rho(T)-\delta_{p,T}](1-e^{K_{2}T})}\right]},
\end{equation}
for $x,y\in H$,$x-y\in\Dom(B_{0}^{-1})$.

\qed

\section{Application}
In this section, we give some simple applications of Theorem \ref{mainthm}.
\begin{corollary}
Assume that $F$, $B$ are determined and independent of t and (H\ref{itemH1}) to (H\ref{itemH5}) hold. If $\lambda_{0}>0$, $\lambda_{0}>K_{1}^{2}+2K_{1}$ and $B(0)\in L_{HS}(H)$, then
\begin{enumerate}[(1)]
\item $P_{t}$ has uniqueness invariant measure $\mu$ and has full support on $H$, $\mu(V)=1$.\label{itemc1}

\item If $\sup_{x}||B(x)||<\infty$, then $\mu(\e^{\epsilon_{0}||\cdot||_{H}^{2}})<\infty$ for some $\epsilon_{0}>0$.\label{itemc2}

\item If there exists $q>0$ such that $\inf_{n}b_{n}^{2q}\lambda_{n}^{q-1}>0$, then $\mu$ has full support on on $H_{0}$.\label{itemc3}
\end{enumerate}
\end{corollary}
\noindent\emph{Proof.}  Let $(V,\ ||\cdot||_{V})=(\Dom(A^{\frac{1}{2}}),\ ||A^{\frac{1}{2}}\cdot||$. Since $\lambda_{0}>0$ and $B(0)\in L_{HS}(H)$, by (H\ref{itemH1}), equation (\ref{equ1}) has strong solution and $P_{t}$ is Feller semigroup. By Ito's formula and $\lambda_{0}>K^{2}_{1}-2K_{1}$, there exists a constant $c>0$ such that
\begin{equation*}
\d ||x_{t}||^{2} \leq \left( c -2(1-\frac{K^{2}_{1}+2K_{1}}{\lambda_{0}})||x_{t}||^{2}_{V} + 2||F(0)||\cdot||x_{t}||\right)\d t + 2\<B(x_{t})\d W_{t},x_{t}\>
\end{equation*}
and
\begin{equation*}
\begin{split}
\d \e^{\epsilon ||x_{t}||^{2}}\leq &\ \epsilon\e^{\epsilon ||x_{t}||^{2}}\left(c-2(1-\frac{K^{2}_{1}+2K_{1}}{\lambda_{0}})||x_{t}||^{2}_{V}
+\frac{\epsilon^{2}}{4}||B^{*}(x_{t})x_{t}||^{2}+2||F(0)||\cdot||x_{t}||\right)\d t \\
&\ + 2\epsilon\e^{\epsilon ||x_{t}||^{2}}\<B(x_{t})\d W_{t},x_{t}\>,
\end{split}
\end{equation*}
for sufficient small $\epsilon$, by H\"{o}lder inequality and noting that $||\cdot||_{V}$ is compact function on $H$, then by standard argument in Theorem 1.2 in \cite{Wang2007}, one can prove (\ref{itemc1}) and (\ref{itemc2}). For (\ref{itemc3}), $\inf_{n}b_{n}^{2q}\lambda_{n}^{q-1}>0$ implies that there exists a constant $c(m)>0$ such that
\begin{equation}
||\cdot||^{2}_{H_{0}} \leq c(m)||\cdot||^{2}+ \frac{1}{m}||\cdot||_{V}^{2},\ \forall m\geq1,
\end{equation}
by Ito's formula, one can get following inequality,
\begin{equation}
\d ||x_{t}(x)-x||^{2} \leq -||x_{t}(x)-x||_{V}^{2}\d t + (c_{1}+c_{2}||x_{t}(x)||^{2})\d t + 2\<B(x_{t})\d W_{t},x_{t}-x\>
\end{equation}
here we denote $x_{t}(x)$ for the process starts from $x$, $c_{1}$, $c_{2}$ are constants depend on $x$. Using Harnack inequality (\ref{Harineq}), (\ref{itemc3}) can be proved following the line of \cite{WangX2011}.
\qed
\begin{corollary}
Assume (H\ref{itemH1}) to (H\ref{itemH5}) hold, $F$ and $B$ are determined and time independent, then
for any $t>0$, $P_{t}$ is $H_{0}$-strong Feller. Let $\mu$ be the $P_{t}$-subinvariant probability with full support on $H_{0}$ as in \cite{RoWang2010}, then the transition density $p_{t}(x,y)$ w.r.t. $\mu$ satisfies
\begin{equation}
||p_{t}(x,\cdot)||_{L^{p}(\mu)}\leq\ \left\{ \int_{H_{0}}{\exp{\left[-\frac{K_{2}\sqrt{q}(\sqrt{q}-1)||x-y||_{H_{0}}^{2}}
{4\delta_{q}[(\sqrt{q}-1)\rho-\delta_{q}](1-e^{K_{2}t})}\right]}\mu(dy)}\right\}^{-\frac{1}{q}}
\end{equation}
for all $1<p<\frac{(K_{3}+\rho)^{2}}{(K_{3}+\rho)^{2}-1}$, here $q=\frac{p}{p-1}$.
\end{corollary}
\noindent\emph{Proof.}  It follows the proof of \cite{Wang2007,RoWang2010,WangX2011}.

\bigskip

\textbf{Acknowledgement}
~The author would like to thank Professor Zdzislaw Brzezniak to provide him the article \cite{Brzez97}, and Professor Feng-Yu Wang for his useful comments.

\paragraph{\Large{Appendix}}
\paragraph{A. Proof of Remark \ref{remark1}}\
\bigskip
\\ \emph{Proof of (1):} since $\bigcup_{n} H_{n}$ is a core of $B_{0}^{-2}$, for any $x\in \mathscr{D}(B_{0}^{-2})$, choose $\{x_{n}\}$ such that $x_{n}\rightarrow x$ and $B_{0}^{-2}x_{n}\rightarrow B_{0}^{-2}x$, hence $B_{0}^{-1}x_{n}\rightarrow x$, as $n \rightarrow +\infty$. Similarly, a sequence $\{y_{n}\}$ with the same property. Therefore
\begin{align*}
&||B_{0}^{-1}[(B(t,x_{n})-B(t,y_{n}))-(B(t,x_{m})-B(t,y_{m})]||_{HS}^{2}\\
\leq &2K_{2}(||B_{0}^{-1}(x_{n}-x_{m})||^{2}
+||B_{0}^{-1}(y_{n}-y_{m})||^{2})
-4\<F(t,x_{n})-F(t,x_{m}),B^{-2}_{0}(x_{n}-x_{m})\>\\
&-4\<F(t,y_{n})-F(t,y_{m}),B^{-2}_{0}(y_{n}-y_{m})\>,
\end{align*}
by the continuous of $F$, we have that $\{B(t,x_{n})-B(t,y_{n})\}$ forms a Cauchy sequence in $L_{HS}(H,H_{0})$. Note that $B(t,x_{n})-B(t,y_{n})$ convergent to $B(t,x)-B(t,y)$ in $L_{HS}(H)$, and $B_{0}^{-1}$ is closed, we have $B(t,x)-B(t,y)\in L_{HS}(H,H_{0})$,
\begin{equation*}
\lim_{n\rightarrow +\infty}(B(t,x_{n})-B(t,y_{n}))=B(t,x)-B(t,y),
\end{equation*}
and
\begin{equation*}
2\<F(t,x)-F(t,y),B_{0}^{-2}(x-y)\>+||B_{0}^{-1}(B(t,x)-B(t,y))||_{HS}^{2}\leq K_{2}||B_{0}^{-1}(x-y)||^{2}.
\end{equation*}
\\\emph{Proof of (2):} we assume $\rho(t)=1$, by definition, it's clear that $B_{0}$ is one to one and has dense range. $$B(t,x)B(t,x)^{*}\geq B_{0}^{2} \Leftrightarrow ||B(t,x)^{*}y||\geq||B_{0}y||, \forall y\in H, $$ implies that $\mathrm{Ran}B(t,x) \supset \mathrm{Ran}B_{0}\mbox{ by Proposition B.1 in \cite{DPZ1992}, and }$  $$||z||\geq||B_{0}(B(t,x)^{*})^{-1}z||, \forall z \in \Ran(B(t,x)^{*}).$$
Since for any $z\in \Ran(B(t,x)^{*})$, $y \in \Ran(B(t,x))$, we have
\begin{equation}
\<B(t,x)^{-1}y,z\>=\<B(t,x)B(t,x)^{-1}y,\ (B(t,x)^{*})^{-1}z\>=\<y,(B(t,x)^{*})^{-1}z\>,
\end{equation}
then
\begin{equation}
z\in\mathscr{D}((B(t,x)^{-1})^{*}),\ (B(t,x)^{-1})^{*}z=(B(t,x)^{*})^{-1}z.
\end{equation}
On the other hand, for any $z\in \mathscr{D}((B(t,x)^{-1})^{*})$, there exists $z^{*}$ such that
\begin{equation}
\<B(t,x)^{-1}y,z\>=\<y,z^{*}\>,\ \forall y \in \mathscr{D}((B(t,x)^{-1})),
\end{equation}
let $u=B(t,x)^{-1}y$, then $\<u,z\>=\<B(t,x)u,z^{*}\>$, we have $z=B(t,x)^{*}z^{*}$ and $$(B(t,x)^{*})^{-1}z=z^{*}=(B(t,x)^{-1})^{*}z,$$ hence $\mathscr{D}((B(t,x)^{-1})^{*})=\mathscr{D}((B(t,x)^{*})^{-1})$. Therefore, $||z||\geq ||B_{0}(B(t,x)^{-1})^{*}z||$, for all $z\in\Ran B(t,x)^{*}$. Since $\Ran(B(t,x)^{*})$ is dense in $H$, $B_{0}(B(t,x)^{-1})^{*}$ can be extended to be a bounded operator on $H$, and for all $z\in H,\ y\in H$, there is $\{z_{n}\}_{n=1}^{+\infty}$, $\lim_{n}z_{n}=z$, such that $\lim_{n}B_{0}(B(t,x)^{-1})^{*}z_{n}=B_{0}(B(t,x)^{-1})^{*}z$, then
\begin{equation}
\begin{split}
&\<B_{0}(B(t,x)^{-1})^{*}z,y\>=\lim_{n}\<B_{0}(B(t,x)^{-1})^{*}z_{n},y\>\\
=&\lim_{n}\<z_{n},(B(t,x)^{-1})B_{0}y\>=\<z,(B(t,x)^{-1})B_{0}y\>,
\end{split}
\end{equation}
hence $||(B(t,x)^{-1})B_{0}y||\leq ||y||$, for all $y\in H$, let $z=B_{0}y$, then $||(B(t,x)^{-1})z||\leq ||B_{0}^{-1}z||$, for all $z\in \Dom(B_{0}^{-1})$. By Proposition B.1 in \cite{DPZ1992}, and the proof above, the converse is easy.

\qed

\paragraph{B. For Lemma \ref{lemma_strongsolution}}\
\bigskip
\\(1) \emph{For local monotonicity}. For any $v_{1},v_{2}\in V$,
\begin{equation}
-2 _{V^{*}}\<A_{0,\epsilon}(v_{1}-v_{2}),v_{2}\>_{V}= -2||\sqrt{A_{0,\epsilon}}(v_{1}-v_{2})||_{H_{0}}^{2}=-2||v_{1}-v_{2}||_{V}^{2},
\end{equation}
\begin{equation}
\begin{split}
&2_{V^{*}}\<F(t,v_{1},v_{2}),v_{1}-v_{2}\>_{V}
+||\hat{B}(t,v_{1})-\hat{B}(t,v_{2})||_{L_{HS}(H,H_{0})}^{2}\\
=&2\<F(t,v_{1},v_{2}),B_{0}^{-2}(v_{1}-v_{2})\>
+||B_{0}^{-1}(\hat{B}(t,v_{1})-\hat{B}(t,v_{2}))||_{HS}^{2}\\
\leq & K_{2}||v_{1}-v_{2}||_{H_{0}}^{2}
\end{split}
\end{equation}
and
\begin{equation}
\begin{split}
&\frac{1}{\xi_{t}}\ _{V^{*}}\<\hat{B}(t,v_{1})G_{n}(t,v_{1})-\hat{B}(t,v_{2})G_{n}(t,v_{2}),v_{1}-v_{2}\>_{V}\\
=&\frac{1}{\xi_{t}}\ \<(\hat{B}(t,v_{1})-\hat{B}(t,v_{2}))G_{n}(t,v_{1})-\hat{B}(t,v_{2})G_{n}(t,v_{1},v_{2}),B_{0}^{-2}(v_{1}-v_{2})\>\\
\leq& \frac{1}{\xi_{t}}||B_{0}^{-1}(\hat{B}(t,v_{1})-\hat{B}(t,v_{2}))||\cdot||G_{n}(t,v_{1})||\cdot||v_{1}-v_{2}||_{H_{0}}\\
&+\frac{1}{\xi_{t}}||B_{0}^{-1}\hat{B}(t,v_{2})||\cdot||G_{n}(t,v_{1},v_{2})||\cdot||v_{1}-v_{2}||_{H_{0}},
\end{split}
\end{equation}
note that, by (H\ref{itemH1}),
\begin{equation*}
\begin{split}
&||B_{0}^{-1}(\hat{B}(t,v_{1})-\hat{B}(t,v_{2}))||_{HS}^{2}\leq K_{2}||v_{1}-v_{2}||_{H_{0}}^{2}-2\<F(t,v_{1},v_{2}),B_{0}^{-2}(v_{1}-v_{2})\>\\
\leq &K_{2}||v_{1}-v_{2}||_{H_{0}}^{2}+
2K_{1}||B_{0}A_{0,\epsilon}^{-\frac{1}{2}}A_{0,\epsilon}^{\frac{1}{2}}(v_{1}-v_{2})||_{H_{0}}
\cdot||B_{0}^{-1}A_{0,\epsilon}^{-\frac{1}{2}}A_{0,\epsilon}^{\frac{1}{2}}(v_{1}-v_{2})||_{H_{0}}\\
\leq &K_{2}||v_{1}-v_{2}||_{H_{0}}^{2}+2K_{1}\left(\sup_{n}\frac{b_{n}}{\sqrt{\lambda_{n}+\epsilon b_{n}^{-2}}}\right)\left(\sup_{n}\frac{1}{b_{n}\sqrt{\lambda_{n}+\epsilon b_{n}^{-2}}}\right)||v_{1}-v_{2}||^{2}_{V}\\
\leq &K_{2}||v_{1}-v_{2}||_{H_{0}}^{2} + \frac{2}{\epsilon}K_{1}||B_{0}||^{2}||v_{1}-v_{2}||_{V}^{2},
\end{split}
\end{equation*}
hence
\begin{equation}
\begin{split}
&\frac{1}{\xi_{t}}||B_{0}^{-1}(\hat{B}(t,v_{1})-\hat{B}(t,v_{2}))||\cdot||G_{n}(t,v_{1})||\cdot||v_{1}-v_{2}||_{H_{0}}\\
\leq & \frac{n}{\xi_{t}}(\sqrt{K_{2}}||v_{1}-v_{2}||_{H_{0}}+\sqrt{\frac{2}{\epsilon}K_{1}}||B_{0}||\cdot||v_{1}-v_{2}||_{V})||v_{1}-v_{2}||_{H_{0}}\\
\leq &(\frac{n}{\xi_{t}}\sqrt{K_{2}}+\frac{n^{2}K_{1}||B_{0}||^{2}}{\epsilon \xi_{t}^{2}\delta^{2}})||v_{1}-v_{2}||_{H_{0}}^{2}+\delta^{2}||v_{1}-v_{2}||_{V}^{2},
\end{split}
\end{equation}
and
\begin{equation}
\begin{split}
&\frac{1}{\xi_{t}}||B_{0}^{-1}\hat{B}(t,v_{2})||\cdot||G_{n}(t,v_{1},v_{2})||\cdot||v_{1}-v_{2}||_{H_{0}}\\
\leq &\frac{1}{\xi_{t}}(\sqrt{K_{2}}||v_{2}||_{H_{0}}
+\sqrt{\frac{2K_{1}}{\epsilon}}||B_{0}||\cdot||v_{2}||_{V})||v_{1}-v_{2}||^{2}_{H_{0}},
\end{split}
\end{equation}
therefore, we have
\begin{equation*}
\begin{split}
&2 _{V^{*}}\<A_{n,\epsilon}(t,v_{1})-A_{n,\epsilon}(t,v_{2}),v_{1}-v_{2}\>_{V}
+||\hat{B}(t,x_{t}-v_{2})-\hat{B}(t,x_{t}-v_{1})||_{L_{HS}(H,H_{0})}^{2}\\
\leq&\left[K_{2}+\frac{2n\sqrt{K_{2}}-2}{\xi_{t}}
+\frac{n^{2}K_{1}||B_{0}||^{2}}{\epsilon^{2}\xi_{t}^{2}\delta^{2}}
+\frac{2}{\xi_{t}}(\sqrt{K_{2}}||v_{2}||_{H_{0}}^{2}
+\sqrt{\frac{2K_{1}}{\epsilon}}||B_{0}||\cdot||v_{2}||_{V}^{2})\right]\times\\
&\times||v_{1}-v_{2}||_{H_{0}}^{2}-2(1-\delta^{2})||v_{1}-v_{2}||_{V}^{2}.
\end{split}
\end{equation*}
\\(2) \emph{For coercivity:}
\begin{equation}
-2 _{V^{*}}\<A_{0,\epsilon}v,v\>_{V}=-2||v||_{V}^{2},\ ||B_{0}^{-1}\hat{B}(t,v)||_{HS}^{2} + 2\<F(t,v,0),B_{0}^{-2}v\>\leq K_{2}||v||^{2},
\end{equation}
\begin{equation}
\begin{split}
\frac{2}{\xi_{t}}\ _{V^{*}}\<\hat{B}(t,v)G_{n}(t,v),v\>_{V}
\leq &\frac{2}{\xi_{t}}||B_{0}^{-1}\hat{B}(t,v)||\cdot||G_{n}(t,v)||\cdot||v||_{H_{0}}\\
\leq&\frac{2n}{\xi_{t}}(K_{2}||v||^{2}+\frac{2K_{1}}{\epsilon}||B_{0}||^{2}||v||_{V}^{2})^{\frac{1}{2}}||v||_{H_{0}}\\
\leq&\frac{2n}{\xi_{t}}(\sqrt{K_{2}}||v||_{H_{0}}+\sqrt{\frac{2K_{1}}{\epsilon}}||B_{0}||\cdot||v||_{V})||v||_{H_{0}}\\
\leq&(\frac{2n\sqrt{K_{2}}}{\xi_{t}}+\frac{2n^{2}K_{1}||B_{0}||^{2}}{\epsilon \xi_{t}^{2}\delta^{2}})||v||_{H_{0}}^{2}+\delta^{2}||v||_{V}^{2},
\end{split}
\end{equation}
hence
\begin{equation}
\begin{split}
&2 _{V^{*}}\<A_{n}(t,v),v\>_{V} + ||\hat{B}(t,v)||_{L_{HS}(H,H_{0})}^{2}\\
\leq& -2(1-\delta^{2})||v||_{V}^{2}+(\frac{n\sqrt{K_{2}}-2}{\xi_{t}}+\frac{n^{2}K_{1}}{\epsilon^{2}\xi_{t}^{2}\delta^{2}})||v||_{H_{0}}^{2}.
\end{split}
\end{equation}
\\(3) \emph{For Growth:}
\begin{equation}
||A_{0,\epsilon}v||^{2}_{V^{*}}=||v||_{V}^{2},\ ||\frac{1}{\xi_{t}}v||_{V^{*}}=\frac{1}{\xi_{t}}||v||_{V^{*}},\ ||F(t,v,0)||_{V^{*}} \leq \frac{K_{1}}{\sqrt{\epsilon}}||v||,
\end{equation}
since, by (H\ref{itemH1}),
\begin{equation}
\begin{split}
| _{V^{*}}\<F(t,v,0),z\>_{V}|=|\<F(t,v,0),B_{0}^{-2}z\>|
\leq K_{1}||v||\cdot||B_{0}^{-2}z||\leq \frac{K_{1}}{\sqrt{\epsilon}}||v||\cdot||z||_{V}.
\end{split}
\end{equation}
And
\begin{equation}
\begin{split}
||\frac{1}{\xi_{t}}\hat{B}(t,v)G_{n}(t,v)||_{V^{*}}
&\leq\frac{||B_{0}||}{\sqrt{\epsilon}\xi_{t}}||\hat{B}(t,v)G_{n}(t,v)||_{H_{0}}\\
&\leq\frac{||B_{0}||}{\sqrt{\epsilon}\xi_{t}}||B_{0}^{-1}\hat{B}(t,v)||\cdot||G_{n}(t,v)||_{L(H_{0},H)}\\
&\leq\frac{||B_{0}||}{\sqrt{\epsilon}\xi_{t}}||(\sqrt{K_{2}}||v||_{H_{0}}
+\sqrt{\frac{2K_{1}}{\epsilon}}||v||_{V}||B_{0}||_{H_{0}})||v||_{H_{0}},
\end{split}
\end{equation}
we have
\begin{equation}
||A_{n,\epsilon}(t,v)||_{V^{*}}^{2}\leq \left(\frac{||B_{0}||^{2}}{\epsilon\xi_{t}}K_{2}
+\left(1+\frac{||B_{0}||^{4}K_{1}}{\epsilon\xi_{t}^{2}}\right)||v||_{V}^{2}\right)(1+||v||_{H_{0}}^{4}).
\end{equation}
\\(4) \emph{For the Lemma 2.2 of \cite{LiuR10}: } We give new estimates to replace inequalities (2.3) and (2.4) there. For convenience, we use the notations there.  In (2.3), we only have to replace $f_{s}\cdot||X_{s}^{(n)}||_{H}^{p-2}$ by $||X_{s}^{(n)}||_{V}||X_{s}^{(n)}||_{H}\cdot||X_{s}^{(n)}||_{H}^{p-2}$ and use the basic inequality
\begin{equation}
a\cdot b \leq \frac{a^{2}}{2\delta}+ \frac{\delta }{2}b^{2}, \forall \delta>0,
\end{equation}
and note that in our case $\alpha = 2$. For (2.4), one can use the following estimate,
\begin{equation*}
\begin{split}
&\E\left(\int_{0}^{\tau_{R}^{n}}{||X_{s}^{(n)}||_{H}^{2p-2}||B(s,X_{s}^{(n)})||_{2}^{2}}\d s\right)^{\frac{1}{2}}\\
\leq&\ \E\left(\int_{0}^{\tau_{R}^{n}}{C||X_{s}^{(n)}
||_{H}^{2p-2}(||X_{s}^{(n)}||_{V}||X_{s}^{(n)}||_{H}+||X_{s}^{(n)}||_{H}^{2})\d s}\right)^{\frac{1}{2}}\\
\leq&\ C(\delta_{1})\E\left(\int_{0}^{\tau_{R}^{n}}{||X_{s}^{(n)}||_{H}^{2p-2}||X_{s}^{(n)}||_{H}^{2}\d s}\right)^{\frac{1}{2}}
+\sqrt{\delta_{1}}\E\left(\int_{0}^{\tau_{R}^{n}}{||X_{s}^{(n)}||_{H}^{2p-2}||X_{s}^{(n)}||_{V}^{2}\d s}\right)^{\frac{1}{2}}\\
\leq&\ \delta_{2}\E\sup_{s\in[0,\tau_{R}^{n}]}{||X_{s}^{(n)}||_{H}^{p}}
+C(\delta_{1},\delta_{2})\E{\int_{0}^{\tau_{R}^{n}}{||X_{s}^{(n)}||_{H}^{p}\d s}}\\
&+\sqrt{\delta_{1}}\E\sup_{s\in[0,\tau_{R}^{n}]}{||X_{s}^{(n)}||_{H}^{\frac{p}{2}}}
\left(\int_{0}^{\tau_{R}^{n}}{||X_{s}^{(n)}||_{H}^{p-2}||X_{s}^{(n)}||_{V}^{2}\d s}\right)^{\frac{1}{2}}\\
\leq&\ (\delta_{2}+\delta_{3})\E\sup_{s\in[0,\tau_{R}^{n}]}{||X_{s}^{(n)}||_{H}^{p}}
+\frac{\delta_{1}}{4\delta_{3}}\E{\int_{0}^{\tau_{R}^{n}}{||X_{s}^{(n)}||_{H}^{p-2}||X_{s}^{(n)}||_{V}^{2}\d s}} +C(\delta_{1},\delta_{2})\E{\int_{0}^{\tau_{R}^{n}}{||X_{s}^{(n)}||_{H}^{p}\d s}},
\end{split}
\end{equation*}
choose $\delta_{2}$, $\delta_{3}$ small enough and $\delta_{1}$ such that $\frac{\delta_{1}}{4\delta_{3}}$ small enough, using $\alpha = 2$ again, then Gronwall's lemma can be applied as in \cite{LiuR10}.

\end{document}